%%%%%%%%%%%%%%%%%%%%%%%%%%%%%%%%%%%%%%%%%%%%%%%%%%%%%%%%%%%%%%%%%%%%%%%%%%%
%%% Title:  A Quantum Group like Structure on Non Commutative 2--Tori
%%% Author: Andreas Cap, Peter W. Michor, Hermann Schichl
%%% Remark: AmSTeX, 5 pages
%%% Series: Letters in Math. Phys. 28 (1993), 251--255
%%%%%%%%%%%%%%%%%%%%%%%%%%%%%%%%%%%%%%%%%%%%%%%%%%%%%%%%%%%%%%%%%%%%%%%%%%%
% TeX-NMB program applied to the file from 1992.11.10;11:45 on 1992.11.10; 11:46
\input amstex
\input amsppt.sty
\hsize 30pc
\vsize 40pc
\magnification=\magstep1
\def\nmb#1#2{#2}         % used for renumbering, TeX should ignore.
\def\totoc{}             %= to table of content, invoked by kms-book.sty
               % for producing index, invoked by kms-book.sty
\def\ign#1{}             %=ignore, invisible entry for the index only

%\pageno251
\redefine\o{\circ}

\define\de{\delta}
\define\ep{\varepsilon}

\define\De{\Delta}

\predefine\ii\i
\redefine\i{^{-1}}
\define\row#1#2#3{#1_{#2},\ldots,#1_{#3}}

%\newsymbol\circledS 1073
\def\today{\ifcase\month\or
 January\or February\or March\or April\or May\or June\or
 July\or August\or September\or October\or November\or December\fi
 \space\number\day, \number\year}
\hyphenation{ho-mo-mor-phism}
\topmatter
\title  
A Quantum Group like Structure on Non Commutative 2--Tori \endtitle
\author  Andreas Cap\\
Peter W. Michor\\
Hermann Schichl
\endauthor
\leftheadtext{\smc Cap, Michor, Schichl}
\rightheadtext{\smc Non commutative two tori}
\thanks{Supported by Project P 7724 PHY 
of `Fonds zur F\"orderung der wissenschaftlichen 
Forschung'\hfill}\endthanks
\affil
Institut f\"ur Mathematik, Universit\"at Wien,\\
Strudlhofgasse 4, A-1090 Wien, Austria.
\endaffil
\address{Institut f\"ur Mathematik, Universit\"at Wien,
Strudlhofgasse 4, A-1090 Wien, Austria.}\endaddress
\email {cap\@pap.univie.ac.at} \endemail
\email {michor\@pap.univie.ac.at} \endemail

%\date{\today}\enddate
%\keywords{}\endkeywords
\subjclass{16W30, 46L87, 58B30, 81R50}\endsubjclass
%\abstract{}\endabstract
\endtopmatter

%\input amspptb.sty
%\userunningheads % used by private style
%\def\leftheadtext{\smc Cap, Michor, Schichl}
%\def\rightheadtext{\smc Non commutative two tori}
%\def\bottremark{\today\hfill}
\document

\head\totoc Introduction \endhead

This paper grew out of a project which aims at a general theory of 
bundles, and in particular principal bundles in non commutative 
differential geometry. The first important question in this direction 
is: What should replace the cartesian product, expressed in terms of 
the algebras representing the factors. Usually one takes the 
(topological) tensor product of the two algebras. The tensor product 
is the categorical coproduct in the category of commutative algebras, 
but does not have a good categorical interpretation in the 
non-commutative case: (coordinates on) the two factors commute; 
why should they?

One can as well consider algebras which have the same underlying 
vector space as the tensor product but are equipped with a 
deformation of the usual multiplication. 

Clearly the question of cartesian products is also relevant for 
studying analogs of Lie groups in non commutative geometry, since the 
comultiplication mappings on the corresponding algebras represent a 
map from the cartesian product of the `group' with itself to the 
`group'. 

In this paper we show that in the case of non commutative two tori 
one gets in a natural way simple structures which have analogous 
formal properties as Hopf algebra structures but with a deformed 
multiplication on the tensor product.

\head\nmb0{1}. Non commutative 2--tori\endhead

The non commutative 2--tori are probably the simplest example of 
non commutative algebras which are commonly thought to describe non 
commutative spaces. They are quite well studied and arise in several 
applications of non commutative geometry to physics, e\.g\. in 
J\.~Bellisard's interpretation of the quantum hall effect. We deal 
here with the smooth version of this algebras, so we consider the 
space $\Cal S(\Bbb Z^2,\Bbb C)$ of all 
complex valued Schwartz sequences (i\.e\. sequences which decay 
faster than any polynomial) $(a_{m,n})$ on $\Bbb Z^2$. Now for a complex 
number $q$ of modulus 1 we define a multiplication on this space as 
follows: Let $\de ^{m,n}$ be the sequence which is one on $(m,n)$ and 
zero on all other points. Then any element of 
$\Cal S(\Bbb Z^2,\Bbb C)$ can be written as $\sum a_{m,n}\de^{m,n}$ 
(this is a convergent sum in the natural Fr\'echet topology on 
$\Cal S(\Bbb Z^2,\Bbb C)$), and we define the multiplication by 
$\de^{m,n}\de^{k,\ell}=q^{-kn}\de^{m+k,n+\ell}$. Now we write 
$U:=\de^{1,0}$ and $V=\de^{0,1}$, then by definition 
$\de^{m,n}=U^mV^n$, and the only relevant relation is $UV=qVU$. Let 
us write $T_q$ for the resulting algebra.

\head\nmb0{2}. The comultiplication 
$\De _q:C^\infty(S^1,\Bbb C)\to T_q$\endhead

First it should be noted that the non commutative 2--tori themselves 
are deformations of the topological tensor product of the algebra 
$C^\infty(S^1,\Bbb C)\cong\Cal S(\Bbb Z,\Bbb C)$ with itself, which 
corresponds to the case $q=1$. In fact the natural comultiplication 
on $\Cal S(\Bbb Z,\Bbb C)$ induced by the group structure on $S^1$ 
given by $(a_n)\mapsto \sum a_n(UV)^n$ makes sense as an algebra 
homomorphism $\De _q:\Cal S(\Bbb Z,\Bbb C)\to T_q$ for any $q$, so 
there is some kind of `multiplication mapping' from any non 
commutative 2--torus to $S^1$. It can be even shown that this 
comultiplication is coassociative in a similar sense as we will use 
below, but there are no counit and antipode mappings fitting to this 
comultiplication.

\heading\totoc\nmb0{3}. The quantum group like structure on $T_q$ 
\endheading

\subhead\nmb.{3.1}\endsubhead
By the general principles of non commutative geometry the algebra 
$T_q$ should be a description of the orbit space $S^1/\Bbb Z$, where 
the action of $\Bbb Z$ on $S^1$ is induced by multiplication by $q$. 
But this space clearly is a group so this should be reflected by some 
kind of coalgebra structure on the algebra $T_q$. 

First we construct the algebra which we consider as a representative 
of the cartesian product of the non commutative torus described by 
$T_q$ with itself. As indicated in the introduction we take as the 
underlying vector space of this algebra the tensor product, i\.e\. 
the space $\Cal S(\Bbb Z^4,\Bbb C)$ of Schwartz sequences on 
$\Bbb Z^4$. As above we write $\de^{k,\ell,m,n}$ for the sequence 
which is one on $(k,\ell,m,n)$ and zero everywhere else. Then we 
define a multiplication on $\Cal S(\Bbb Z^4,\Bbb C)$ by:
$$
\de^{k_1,\ell_1,m_1,n_1}\de^{k_2,\ell_2,m_2,n_2}=
q^{\tfrac{k_2n_1}2-k_2\ell_1-\tfrac{m_1\ell_2}2-m_1n_2}
\de^{k_1+k_2,\ell_1+\ell_2,m_1+m_2,n_1+n_2}
$$
Writing $U_1:=\de^{1,0,0,0}$, $V_1:=\de^{0,1,0,0}$, 
$U_2:=\de^{0,0,1,0}$ and $V_2:=\de^{0,0,0,1}$ we see that 
$\de^{k,\ell,m,n}=U_1^kV_1^\ell U_2^mV_2^n$ and the relations between 
these generators are:
$$\alignat2
U_1V_1&=qV_1U_1 & U_2V_2&=qV_2U_2\\
U_1V_2&=q^{-1/2}V_2U_1 \quad & \quad V_1U_2&=q^{1/2}U_2V_1\\
U_1U_2&=U_2U_1 & V_1V_2&=V_2V_1
\endalignat$$
Let us remark a little on these relations. Clearly we should have two 
canonical copies of the original algebra $T_q$ contained, so the 
relation between $U_1$ and $V_1$ as well as the one between $U_2$ and 
$V_2$ is clear. Then it turns out that all other relations fall out 
of the further development if one assumes that any product of two of 
the generators is a scalar multiple of the opposite product. We write 
$P_q^2$ for the resulting algebra in the sequel.

\subhead\nmb.{3.2}\endsubhead
Now we define the comultiplication $\De :T_q\to P_q^2$ by 
$\De (U)=U_1U_2$ and $\De (V)=V_1V_2$. Then this induces an algebra 
homomorphism, which is obviously continuous, since 
$$U_1U_2V_1V_2=q^{-1/2}U_1V_1U_2V_2=q^{3/2}V_1U_1V_2U_2=qV_1V_2U_1U_2.$$
To formulate coassociativity we need a representative of the three 
fold product of a non commutative torus with itself. Again we want 
this to be a deformation of the third tensor power of $T_q$, so the 
underlying vector space is $\Cal S(\Bbb Z^6,\Bbb C)$. We write down 
this algebra only in terms of generators and relations, so we write 
a Schwartz sequence as 
$\sum a_{i,j,k,\ell ,m,n}U_1^iV_1^jU_2^kV_2^\ell U_3^mV_3^n$. The 
choice of relations in this algebra is dictated by requiring the 
same deformation which led to $P^2_q$ in the first and second factor 
as well as in the second and third factor and no additional 
deformations. This gives the following relations:
$$\alignat3
U_1V_1&=qV_1U_1 & U_2V_2&=qV_2U_2 & U_3V_3&=qV_3U_3 \\
U_1V_2&=q^{-1/2}V_2U_1 \quad &\quad U_2V_3&=q^{-1/2}V_3U_2 
\quad &\quad U_3V_1&=V_1U_3 \\
V_1U_2&=q^{1/2}U_2V_1 & V_2U_3&=q^{1/2}U_3V_2 & V_3U_1&=U_1V_3 \\
U_1U_2&=U_2U_1 & U_2U_3&=U_3U_2 & U_3U_1&=U_1U_3 \\
V_1V_2&=V_2V_1 & V_2V_3&=V_3V_2 & V_3V_1&=V_1V_3 
\endalignat$$
We write $P^3_q$ for the resulting algebra.

\subhead\nmb.{3.3}\endsubhead
Taking into account the canonical vector space isomorphisms 
$P^2_q\cong T_q\hat\otimes T_q$ and 
$P^3_q\cong T_q\hat\otimes T_q\hat\otimes T_q$, where $\hat\otimes$ 
denotes the projective tensor product, we get continuous linear maps 
$(\De,Id)$ and $(Id,\De)$ from $P^2_q$ to $P^3_q$ which are induced 
by $\De\hat\otimes Id$ and $Id\hat\otimes \De$, respectively. Now it 
turns out that in this setting these maps are even algebra 
homomorphisms. We show this only for $(\De,Id)$, the proof for 
$(Id,\De)$ is completely analogous.

On the generators of $P^2_q$ the map $(\De,Id)$ is given by
$$U_1\mapsto U_1U_2\quad V_1\mapsto V_1V_2\quad U_2\mapsto U_3\quad 
V_2\mapsto V_3,$$
and this induces an algebra homomorphism since
$$\gather
U_1U_2V_3=q^{-1/2}U_1V_3U_2=q^{-1/2}V_3U_1U_2\\
V_1V_2U_3=q^{1/2}V_1U_3V_2=q^{1/2}U_3V_1V_2\\
U_1U_2U_3=U_3U_1U_2\\
V_1V_2V_3=V_3V_1V_2
\endgather$$
But then coassociativity is obvious since both $(\De,Id)\o\De$ and 
$(Id,\De)\o\De$ are continuous algebra homomorphisms $T_q\to P^3_q$ which 
map $U$ to $U_1U_2U_3$ and $V$ to $V_1V_2V_3$.

\subhead\nmb.{3.4}\endsubhead
Let us next turn to the counit $\ep :T_q\to \Bbb C$. We cannot expect 
to get a counit which is an algebra homomorphism since there are no 
nonzero homomorphisms from $T_q$ to a commutative algebra. This can be 
interpreted as the fact that the non commutative torus has no 
classical points. So we define $\ep$ as a linear map by 
$\ep (U^kV^\ell):=q^{\tfrac{k\ell}2}$. Clearly this defines a 
continuous linear mapping. Again using the canonical vector space 
isomorphism $P^2_q\cong T_q\hat\otimes T_q$ we get continuous linear 
mappings $(\ep,Id)$ and $(Id,\ep)$ from $P^2_q$ to $T_q$. These maps 
are given by 
$(\ep,Id)(U_1^kV_1^\ell U_2^mV_2^n)=q^{\tfrac{k\ell}2}U^mV^n$ and by 
$(Id,\ep)(U_1^kV_1^\ell U_2^mV_2^n)=q^{\tfrac{mn}2}U^kV^\ell$.
To prove that $\ep$ is in fact a counit for the comultiplication 
$\De$ we have to show that both $(\ep,Id)\o\De$ and $(Id,\ep)\o\De$ 
are the identity map. Since these are continuous linear maps it 
suffices to check this on elements of the form $U^kV^\ell$. But
$$\De(U^kV^\ell)=(U_1U_2)^k(V_1V_2)^\ell=U_1^kU_2^kV_1^\ell V_2^\ell =
q^{-\tfrac{k\ell}2}U_1^kV_1^\ell U_2^kV_2^\ell,$$ so the result is 
obvious.

\subhead\nmb.{3.5}\endsubhead
Finally we define the antipode map $S:T_q\to T_q$. The obvious choice 
is $S(U):=U^{-1}$ and $S(V):=V^{-1}$. Then this extends to a 
continuous algebra homomorphism since $U^{-1}V^{-1}=qV^{-1}U^{-1}$. 
So in contrast to the case of Hopf algebras we get an antipode which 
is a homomorphism and not an anti homomorphism. This seems rather 
positive since algebra homomorphisms should correspond to mappings of 
the underlying `spaces' while the meaning of anti homomorphisms is 
rather unclear. 

As before we get continuous linear mappings $(S,Id)$  and $(Id,S)$ 
from $P^2_q$ to itself using the vector space isomorphism  with the 
tensor product. These maps are given by 
$U_1^kV_1^\ell U_2^mV_2^n\mapsto U_1^{-k}V_1^{-\ell}U_2^mV_2^n$ and
$U_1^kV_1^\ell U_2^mV_2^n\mapsto U_1^kV_1^{\ell}U_2^{-m}V_2^{-n}$, 
respectively, and they are not algebra homomorphisms (since the 
multiplication on $P^2_q$ is twisted). Moreover the same method leads 
to a continuous linear map $\mu :P^2_q\to T_q$ which represents the 
multiplication of $T_q$ and is given by 
$\mu (U_1^kV_1^\ell U_2^mV_2^n)=U^kV^\ell U^mV^n=
q^{-\ell m}U^{k+m}V^{\ell +n}$. 

To prove that the antipode is really an analog of a group inversion 
we have to show that $(\mu\o(S,Id)\o\De)(x)=\ep (x)\cdot 1$ and 
likewise with $(S,Id)$ replaced by $(Id,S)$. Since all these maps are 
linear and continuous it suffices to check this on elements of the 
form $U^kV^\ell$. But for these we get:
$$\align
(\mu\o(S,Id)\o\De)(U^kV^\ell)&=
	(\mu\o(S,Id))(q^{-\tfrac{k\ell}2}U_1^kV_1^\ell U_2^kV_2^\ell)\\
&= \mu(q^{-\tfrac{k\ell}2}U_1^{-k}V_1^{-\ell}U_2^kV_2^\ell)\\
&=q^{-\tfrac{k\ell}2}q^{k\ell}=q^{\tfrac{k\ell}2}=\ep(U^kV^\ell)\cdot 1
\endalign$$
and likewise with $(Id,S)$.

\enddocument